\documentclass[12pt]{article}
\usepackage{amsfonts}
\newtheorem{theorem}{Theorem}

\newtheorem{definition}[theorem]{Definition}
\newtheorem{example}[theorem]{Example}

\newtheorem{remark}[theorem]{Remark}
\newenvironment{proof}[1][Proof]{\textbf{#1.} }{\ \rule{0.5em}{0.5em}}

\title{An Introduction to Some Spaces of Interval Functions}
\author{Roumen Anguelov\\
Department of Mathematics and Applied Mathematics\\
University of Pretoria\\
SOUTH AFRICA\\
anguelov@scientia.up.ac.za\vspace{3mm}}
\date{}

\begin{document}
\maketitle
 \begin{abstract}
The paper gives a brief account of the spaces of interval
functions defined through the concepts of H-continuity,
D-continuity and S-continuity. All three continuity concepts
generalize the usual concept of continuity for real (point valued)
functions. The properties of the functions in these new spaces are
discussed and investigated, preserving essential properties of the
usual continuous real functions being of primary interest. Various
ways in which the spaces of H-continuous, D-continuous and
S-continuous interval functions complement the spaces of
continuous real functions are discussed.
 \end{abstract}


\section{Introduction}

The aim of this paper is to give an introduction to the spaces of
interval functions which have emerged in connection with
applications to real analysis, approximation theory and partial
differential equations. These spaces are all based on extending
the concept of continuity of real functions to interval functions.
We denote by $\mathbb{A}(\Omega )$ the set of all functions
defined on an open set $\Omega\subset\mathbb{R}^{n}$ with values
which are finite or infinite closed real intervals, that is,
\[
\mathbb{A}(\Omega )=\{f:\Omega \rightarrow
\mathbb{I\,\overline{\mathbb{R}} }\},
\]%
where
$\mathbb{I\,}\overline{\mathbb{R}}=\{[\underline{a},\overline
{a}]:\underline{a},\overline{a}\in\overline{\mathbb{R}}=\mathbb{R\cup
\{\pm\infty\}},$ $\underline{a}\leq\overline{a}\}$. Given an
interval
$a=[\underline{a},\overline{a}]\in\mathbb{I\,\overline{\mathbb{R}}}$,
\[
w(a)=\left\{
\begin{tabular}
[c]{lll}%
$\overline{a}-\underline{a}$ & if & $\underline{a},\overline{a}$ finite,\\
$\infty$ & if & $\underline{a}<\overline{a}=\infty$ or $\underline{a}%
=-\infty<\overline{a}$,\\
0 & if & $\underline{a}=\overline{a}=\pm\infty,$%
\end{tabular}
\right.
\]
is the width of $a$. An extended real interval $a$ is called a
proper interval if $w(a)>0$ and degenerate or point interval if
$w(a)=0$. Identifying $a\in \overline{\mathbb{R}}$ with the
degenerate interval $[a,a]\in \mathbb{I}\,\overline{\mathbb{R}}$,
we consider $\overline{\mathbb{R}}$ as a subset of
$\mathbb{I}\,\overline{\mathbb{R}}$. In this way
$\mathbb{A}(\Omega )$ contains the set of extended real valued functions, namely,%
\[
\mathcal{A}(\Omega)=\{f:\Omega \rightarrow
\overline{\mathbb{R}}\}.
\]
The set of all continuous real point valued functions $C(\Omega)$
is a subset of $\mathcal{A}(\Omega)$, that is, we have the
inclusions
\[C(\Omega)\subseteq \mathcal{A}(\Omega) \subseteq
\mathbb{A}(\Omega).\]

For every $x\in \Omega$ the value of $f\in \mathbb{A}(\Omega)$ is
an interval $[\underline{f}(x),\overline{f}(x)]$. Hence, the
function $f$ can be written
in the form $f=[\underline{f},\overline{f}]$, where $\underline{f},\overline{f%
}\in \mathcal{A}(\Omega)$ and $\underline{f}(x)\leq
\overline{f}(x),\;x\in\Omega$. A trivial and direct extension of
the concept of continuity of real functions to interval functions
is given in the next definition.
\begin{definition}\label{Defcont}
A function $f=[\underline{f},\overline{f}]\in\mathbb{A}(\Omega )$
is called continuous if the functions $\underline{f}$ and
$\overline{f}$ are both continuous extended real valued functions.
\end{definition}
If the values of $f$ are finite intervals then
\[
f-{\rm
continuous}\;\Longleftrightarrow\;\underline{f},\overline{f}\in
C(\Omega).
\]

An important issue in interval analysis is the construction of
enclosures. Let ${\cal F}$ be a set of continuous real function on
$\Omega$, i.e., ${\cal F}\subset C(\Omega)$. A function
$\psi\in\mathbb{A}(\Omega )$ is called a continuous interval
enclosure of ${\cal F}$ if $\psi$ is continuous and
\[
\phi(x)\in \psi(x),\;x\in\Omega,\;\phi\in {\cal F}.
\]
A continuous interval enclosure $\psi$ of {\cal F} is called
minimal if for any continuous interval function $h$ we have
\begin{equation}\label{hullminprop}
\phi(x)\in h(x),\;x\in\Omega,\;\phi\in{\cal F} \Longrightarrow
\psi(x)\subseteq h(x),\;x\in\Omega.
\end{equation}
Clearly if the minimal continuous interval enclosure exists it is
given by the following function also called interval hull of
${\cal F}$:
\begin{equation}\label{hull}%
{\rm hull}({\cal F})(x)=\bigcap_{\psi\in \hat{\cal
F}}\psi(x),\;x\in\Omega,
\end{equation}
where $\hat{\cal F}$ is the set of all continuous interval
enclosures of ${\cal F}$. The following example shows that a
continuous minimal interval enclosure of a set of continuous
functions does not always exist, that is, the interval hull of a
set of continuous functions is not necessarily a continuous
functions.
\begin{example}\label{ExHullNotCont}
Let us consider the following subset of $C(\mathbb{R})$
\[{\cal F}=\{\phi_\lambda:\lambda>0\}
\]
where
\[\phi_\lambda(x)=\left\{\begin{tabular}{ccc}
$\displaystyle\frac{1-e^{-\lambda x}}{1+e^{-\lambda x}}$&if&$x\geq
0$\\
$0$&if&$x<0$
\end{tabular}\right.
\]
The interval hull of ${\cal F}$
\begin{eqnarray*}
f(x)&=&{\rm hull}({\cal F})(x)= \left\{\begin{tabular}{ccc}
$[0,1]$&if&$x\geq 0$\\ $0$&if&$x< 0$\end{tabular}\right.
\end{eqnarray*}
is not a continuous interval function indicating that the set
${\cal F}$ does not have a minimal continuous interval enclosure.
\end{example}
The fact that the sets of continuous functions do not always have
minimal continuous interval enclosures is not surprising as it
relates to the fact that the set of continuous real functions
$C(\Omega)$ is not Dedeking order complete, that is, the supremum
and the infimum of a bounded set of continuous functions do not
always exist in $C(\Omega)$. The order on $C(\Omega)$ considered
in this regard is the usual one induced by the order on
$\mathbb{R}$ in a point-wise way. In the same way partial order
for interval functions can be point-wise induced by a partial
order on  $\mathbb{I\,\overline{\mathbb{R}}}$.  Note that a
partial order which extends the total order on
$\overline{\mathbb{R}}$ can be defined on
$\mathbb{I\,}\overline{\mathbb{R}}$ in more than one way, see
Appendix 2. However, it proves useful to consider on
$\mathbb{I\,}\overline{\mathbb{R}}$ the partial order $\leq$ defined by%
\begin{equation}
\lbrack\underline{a},\overline{a}]\leq\lbrack\underline{b},\overline
{b}]\Longleftrightarrow\underline{a}\leq\underline{b},\;\overline{a}%
\leq\overline{b}. \label{iorder}%
\end{equation}
The partial order induced in $\mathbb{A}(\Omega)$ by
(\ref{iorder}) in a point-wise
way, i.e.,%
\begin{equation}
f\leq g\Longleftrightarrow f(x)\leq g(x),\;x\in \Omega, \label{forder}%
\end{equation}
is an extension of the usual point-wise order on the set of
extended real valued functions $\mathcal{A}(\Omega)$ and in
particular on the set $C(\Omega)$. A partial order naturally
related to interval spaces is the relation inclusion. Similarly to
(\ref{forder}) it can be defined for interval functions in a
point-wise way, that is,
\begin{equation}
f\subseteq g\Longleftrightarrow f(x)\subseteq g(x),\;x\in \Omega, \label{definclusion}%
\end{equation}
We will see in the sequel that suitable subsets of the set of all
Hausdorff continuous (H-continuous) interval functions considered
with the order (\ref{forder}) are Dedeking order completions of
well known sets of continuous function, thus improving an earlier
result by Dilworth \cite{Dilworth}. We will also define the set of
Dilworth continuous (D-continuous) interval functions which
contains the set of Hausdorff continuous functions as well as all
interval hulls of bounded sets of continuous functions. This set
provides a Dedekind order completion of the set of continuous
interval functions both with respect to the order relation $\leq$
given in (\ref{forder}) and the order relation inclusion
($\subseteq$) given in (\ref{definclusion}). The even larger set
of S-continuous functions contains also the topological completion
of $C(\Omega)$ with respect to the Hausdorff distance between
functions as defined in \cite{Sendov}. The three concepts of
continuity, namely, H-continuity, D-continuity, S-continuity,
considered in this paper in addition to the continuity concept in
Definition \ref{Defcont} all generalize the concept of continuity
of real functions in the sense that if a real function $f$ is
H-continuous, or D-continuous, or S-continuous, then it is
continuous. These concepts are defined through the lower and upper
Baire operators and the graph completion operator. Since these
operators are fundamental for the exposition we recall in the next
section, that is, Section 2, their definitions and essential
properties. The use of extended real intervals is partially
motivated by the fact that the definitions of the Baire operators
involve infimums and supremums which might not exists in the realm
of the usual (finite) real intervals. However, most applications
involve functions which are finite or nearly finite in the sense
of the following definitions.
\begin{definition}\label{Deffinite}
A function $f\in\mathbb{A}(\Omega)$ is called finite if
\[
|f(x)|<\infty,\;x\in\Omega
\]
\end{definition}
\begin{definition}\label{Defnearlyfinite}
A function $f\in\mathbb{A}(\Omega)$ is called nearly finite if
there exists an open and dense subset $D$ of $\Omega$ such that
\[
|f(x)|<\infty,\;x\in D
\]
\end{definition}

The concepts of continuity for interval functions mentioned above
as well as the respective spaces are discussed in Section 3.  The
spaces of Hausdorff continuous interval functions are particularly
considered in Section 4 while Sections 5 and 6 deal with the space
of D-continuous functions and the space of S-continuous functions,
respectively. Various ways in which the considered spaces of
interval functions complement the spaces of continuous functions
are discussed in the Sections 7-9.

\section{The Baire operators and the graph completion operator}

For every $x\in \Omega$, $B_{\delta }(x)$ denotes the open $\delta $%
-neighborhood of $x$ in $\Omega $, that is,
\[
B_{\delta }(x)=\{y\in \Omega :||x-y||<\delta \}.
\]%
Let $D$ be a dense subset of $\Omega$.  The pair of mappings
$I(D,\Omega,\cdot),$ $S(D,\Omega,\cdot):\mathbb{A}(D )\rightarrow
\mathcal{A}(\Omega )$ defined by
\begin{eqnarray}
I(D,\Omega ,f)(x) &=&\sup_{\delta >0}\inf \{z\in f(y):y\in
B_{\delta
}(x)\cap D\},x\in \Omega ,  \label{lbfgen} \\
S(D,\Omega ,f)(x) &=&\inf_{\delta >0}\sup \{z\in f(y):y\in
B_{\delta }(x)\cap D\},x\in \Omega ,  \label{ubfgen}
\end{eqnarray}%
are called lower Baire and upper Baire operators, respectively.
Clearly for every $f\in \mathbb{A}(D)$ we have
\[
I(D,\Omega ,f)(x)\leq f(x)\leq S(D,\Omega ,f)(x),\;x\in\Omega.
\]
Hence the mapping $F:\mathbb{A}(D)\rightarrow \mathbb{A}(\Omega
)$, called a graph completion operator, where
\begin{equation}\label{gcgen}
F(D,\Omega ,f)(x)=[I(D,\Omega ,f)(x),S(D,\Omega ,f)(x)],\;x\in
\Omega ,\;f\in \mathbb{A}(\Omega ),
\end{equation}
is well defined and we have the inclusion
\begin{equation}\label{fFincl}
f(x)\subseteq F(D,\Omega,f)(x),\ x\in\Omega.
\end{equation}
The name of this operator is derived from the fact that
considering the graphs of $f$ and $F(D,\Omega,f)$ as subsets of
the topological space $\Omega\times \overline{\mathbb{R}}$, the
graph of $F(D,\Omega,f)$ is the minimal closed set which is a
graph of interval function on $\Omega$ and contains the the graph
of $f$, see \cite{Sendov}. In the case when $D=\Omega$ the sets
$D$ and $\Omega$ will be usually omitted from the operators'
argument lists, that is,
\[
I(f)=I(\Omega,\Omega,f)\ ,\ \ S(f)=S(\Omega,\Omega,f)\ ,\ \
F(f)=F(\Omega,\Omega,f)
\]

The lower Baire operator, the upper Baire operator and the graph
completion operator are closely connected with the order relations
in the domains of their arguments and their ranges. The following
monotonicity properties with respect to the order relations
(\ref{forder}) and (\ref{definclusion}) follow immediately from
the definitions of the operators.

\begin{itemize}
\item[1.]The lower Baire operator, the upper Baire operator and
the graph completion operator are all monotone increasing with
respect to the their functional argument, that is, if $D$ is a
dense subset of $\Omega$, for every two functions
$f,g\in\mathbb{A}(D)$ we have%
\begin{equation}\label{mon}
f(x)\leq g(x),~x\in D~\Longrightarrow ~ \left\{\begin{tabular}{l}
$I(D,\Omega,f)(x)\leq I(D,\Omega,g)(x),~x\in\Omega$\\
$S(D,\Omega,f)(x)\leq S(D,\Omega,g)(x),~x\in\Omega$\\
$F(D,\Omega,f)(x)\leq F(D,\Omega,g)(x),~x\in\Omega$
\end{tabular}\right.
\end{equation}

\item[2.]The graph completion operator is inclusion isotone with
respect to the functional argument, that is, if $f,g\in
\mathbb{A}(D)$, where $D$ is dense in $\Omega$, then
\begin{equation}\label{Finclmon1}
f(x)\subseteq g(x),\;x\in D\Longrightarrow F(D,\Omega
,f)(x)\subseteq F(D,\Omega ,g)(x),\;x\in \Omega .
\end{equation}

\item[3.]The graph completion operator is inclusion isotone with
respect to the set $D$ in the sense that if $D_{1}$ and $D_{2}$
are dense subsets of $\Omega $ and $f\in \mathbb{A}(D_{1}\cup
D_{2})$ then
\begin{equation}\label{Finclmon2}
D_{1}\subseteq D_{2}\Longrightarrow F(D_{1},\Omega ,f)(x)\subseteq
F(D_{2},\Omega ,f)(x),x\in \Omega .
\end{equation}
This, in particular, means that for any dense subset $D$ of
$\Omega $ and $f\in
\mathbb{A}(\Omega )$ we have%
\begin{equation}\label{Finclmon3}
F(D,\Omega ,f)(x)\subseteq F(f)(x),x\in \Omega .
\end{equation}
\end{itemize}

Further details about the Baire operators and the graph completion
operator are given in Appendix 1 in connection with the
semi-continuous function. From the properties presented there it
can be easily seen that all three operators are idempotent and it
fact we have the following stronger property: If the sets $D_1$
and $D_2$ are both dense in $\Omega$ and $D_1\subseteq D_2$ then
\begin{eqnarray}
I(D_2,\Omega,\cdot)\circ I(D_1,\Omega,\cdot)&=&I(D_1,\Omega,\cdot)\nonumber\\
S(D_2,\Omega,\cdot)\circ S(D_1,\Omega,\cdot)&=&S(D_1,\Omega,\cdot)\label{idemp}\\
F(D_2,\Omega,\cdot)\circ
F(D_1,\Omega,\cdot)&=&F(D_1,\Omega,\cdot)\nonumber
\end{eqnarray}

\section{Three concepts of continuity}

Using the graph completion operator we define the following three
concepts.

\begin{definition}\label{Defscont}
A function $f\in \mathbb{A}(\Omega )$ is called S-continuous, if
$F(f)=f$.
\end{definition}

\begin{definition}\label{Defnormal}
A function $f\in \mathbb{A}(\Omega )$ is called Dilworth
continuous or shortly D-continuous if for every
dense subset $D$ of $\Omega $ we have%
\[
F(D,\Omega ,f)=f .
\]
\end{definition}

\begin{definition}\label{DefHcont}%
A function $f\in\mathbb{A}(\Omega)$ is called Hausdorff
continuous, or H-continuous, if for every function
$g\in\mathbb{A}(\Omega)$ which satisfies the inclusion
$g(x)\subseteq f(x),$ $x\in \Omega$, we have $F(g)(x)=f(x),$ $x\in
\Omega $.
\end{definition}

The following theorem indicates that the H-continuity is the
strongest concept among the three while the S-continuity is the
weakest.
\begin{theorem}\label{thcontregscont}
Let $f\in\mathbb{A}(\Omega)$. Then the following implications hold
\begin{equation}\label{hcontregscont}
f\mbox{ is H-continuous}\Longrightarrow f\mbox{ is
D-continuous}\Longrightarrow f\mbox{ is S-continuous}
\end{equation}
\end{theorem}
\begin{proof}
Let $f$ be H-continuous. From the inclusion \[f(x)\subseteq
f(x),\;x\in\Omega,\] and the Definition \ref{DefHcont} it follows
that
\[F(f)=f\ ,\]
which means that $f$ is S-continuous. Let $D$ be any dense subset
of $\Omega$. Due to the property (\ref{Finclmon3}) of the graph
completion operator we have
\[F(D,\Omega,f)(x)\subseteq F(f)(x)=f(x),\;x\in\Omega
\]
Using the idempotence of the graph completion operator, see
(\ref{idemp}), and the minimality property in the definition of
Hausdorff continuity we have
\[
F(D,\Omega,f)(x)=F(F(D,\Omega,f))=f(x),\;x\in\Omega,
\]
which shows that function $f$ is D-continuous.

The second implication in the theorem follows immediately from
Definition \ref{Defnormal} by taking $D=\Omega$.
\end{proof}

We will use the notations:
\begin{itemize}
\item[] $\mathbb{F}(\Omega)$ - the set of all S-continuous
interval functions defined on $\Omega$,
\item[]$\mathbb{G}(\Omega)$ - the set of all D-continuous interval
functions defined on $\Omega$, \item[] $\mathbb{H}(\Omega)$ - the
set of all H-continuous interval functions defined on $\Omega$,
\end{itemize}
where the term interval functions means extended real interval
valued functions, that is, the elements of $\mathbb{A}(\Omega)$.
The following inclusions follow from (\ref{hcontregscont})
\[\mathbb{H}(\Omega)\subseteq \mathbb{G}(\Omega)\subseteq
\mathbb{F}(\Omega).
\]

With every interval function $f$ one can associate S-continuous,
D-continuous and H-continuous functions as stated in the next
theorem.

\begin{theorem}
\label{tFSI} Let $f\in \mathbb{A}(\Omega)$. Then
\begin{itemize}
\item[(i)] for every dense subset $D$ of $\Omega$ the function $F(D,\Omega,f)$ is S-continuous;%
\item[(ii)] the function $G(f)=[I(S(I(f))),S(I(S(f)))]$ is D-continuous;%
\item[(iii)] both functions $F(S(I(f)))$ and $F(I(S(f)))$ are
H-continuous and
\[
F(S(I(f)))\leq F(I(S(f)))\ .
\]
\end{itemize}
\end{theorem}
\begin{proof}
(i) Follows immediately from the idempotence of the operator $F$,
see (\ref{idemp}).

(iii) Denote $g=F(S(I(f)))$. The function $g$ can also be written
in the form $g=[\underline{g},\overline{g}]$, where
\begin{eqnarray}
\underline{g}  &  =&I(S(I(f))),\label{gisif}\\
\overline{g}  &  =&S(S(I(f)))=S(I(f)). \label{gsif}%
\end{eqnarray}
Let us assume that the function $h=[\underline{h},\overline{h}]
\in\mathbb{A}(\Omega)$ satisfies the inclusion
\[
h(x)\subseteq g(x),\ x\in\Omega,
\]
or, equivalently,
\begin{equation}\label{iiiineq}
I(S(I(f)))\leq \underline{h}\leq \overline{h}\leq S(I(f)).
\end{equation}
Using the idempotence, see (\ref{idemp}), and the monotonicity,
see (\ref{mon}), of the operators $I$ and $S$ we have
\[
\underline{g}=I(S(I(f)))=I(I(S(I(f))))\leq I(\underline{h})\leq
I(\overline{h})\leq I(S(I(f)))=\underline{g}
\]
Therefore
\begin{equation}\label{iiieq1}
\underline{g}=I(\underline{h})
\end{equation}%
In a similar way using the same idempotence and monotonicity
properties with the inequalities (\ref{iiiineq}) we obtain
\begin{eqnarray*}
\overline{g}&=&S(I(f))\ =\ S(S(I(f)))\ \geq\ S(\overline{h})\
\geq\ S(\underline{h})\\
&\geq& S(I(S(I(f))))\ \geq\ S(I(I(I(f))))\ =\ S(I(f))\ =\
\overline{g}
\end{eqnarray*}
Hence
\begin{equation}\label{iiieq2}
\overline{g}=S(\overline{h})
\end{equation}
From (\ref{iiieq1}) and (\ref{iiieq2}) it follows that
\[
g=F(h).
\]
Hence, the function $g=F(S(I(f))$ is H-continuous. The
H-continuity of $F(I(S(f)))$ is proved in the same way. The
inequality in the theorem follows also from the properties
(\ref{idemp}) and (\ref{mon}). We have
\begin{eqnarray*}
F(S(I(f)))&\leq& F(S(I(S(f))))~=~[I(S(I(S(f)))),S(I(S(f)))]\\
&\leq&
[I(S(S(S(f)))),S(I(S(f)))]~=~[I(S(f)),S(I(S(f)))]~=~F(I(S(f)))
\end{eqnarray*}

(ii) Let $f\in\mathbb{A}(\Omega)$ and let $D$ be any dense subset
of $\Omega$. It follows from (iii), which is proved above, that
the function
\[
g=[\underline{g},\overline{g}]=F(S(I(f)))=[I(S(I(f))),S(I(f))]
\]
is H-continuous, and therefore D-continuous, see Theorem
\ref{thcontregscont}. Hence
\[F(D,\Omega,g)=g
\]
which in particular means that
\[I(D,\Omega,\underline{g})=\underline{g}.\]
Therefore
\[I(D,\Omega,I(S(I(f))))=I(S(I(f))).\]
In the same way we prove that
\[S(D,\Omega,S(I(S(f))))=S(I(S(f))).\]
Thus we have
\begin{eqnarray*}
F(D,\Omega,G(f))&=&[I(D,\Omega,I(S(I(f)))),S(D,\Omega,S(I(S(f))))]\\
&=&[I(S(I(f))),S(I(S(f)))]\ =\ G(f),
\end{eqnarray*}
which shows that the function $G(f)$ is D-continuous.
\end{proof}

The operator $G:\mathbb{A}(\Omega)\longrightarrow
\mathbb{G}(\Omega)$ is a useful tool for studying the D-continuous
interval functions. Similar to the graph completion operator it
can be defined for functions in $\mathbb{A}(D)$ where $D$ is a
dense subset of $\Omega$. However, since this will not be used in
this paper we will consider $G$ only on $\mathbb{A}(\Omega)$ as
stated in the following definition.

\begin{definition}\label{Defnorm}
The operator $G:\mathbb{A}(\Omega)\longrightarrow
\mathbb{G}(\Omega)$ given by \[ G(f)=[I(S(I(f))),S(I(S(f)))]\] is
called normalizing operator.
\end{definition}

Theorem \ref{tFSI} is illustrated by the following example.

\begin{example}\label{ExtFSI}
Consider the function $f\in\mathbb{A}(\mathbb{R})$ given by
\[
f(x)=\left\{\begin{tabular}{ccl} $[-1,1]$&if&$x\in\mathbb{Z}$\\
$0$&if&$x\in (-\infty,0)\backslash\mathbb{Z}$\\
$[0,1]$&if&$x\in (0,\infty)\backslash\mathbb{Z}$ \end{tabular}
\right.
\]

We have $F(f)=f$ meaning that $f$ is S-continuous.

The D-continuous function $G(f)$ is given by
\[
G(f)(x)=\left\{\begin{tabular}{ccl}
$0$&if&$x\in (-\infty,0)$\\
$[0,1]$&if&$x\in [0,\infty)$ \end{tabular} \right.
\]

Finally, we have the H-continuous functions
\begin{eqnarray*}
F(S(I(f)))(x)&=&0,\;x\in\mathbb{R}\\
F(I(S(f)))(x)&=&\left\{\begin{tabular}{ccl}
$0$&if&$x\in (-\infty,0)$\\
$[0,1]$&if&$x=0$\\
$1$&if&$x\in (0,\infty)$ \end{tabular} \right.
\end{eqnarray*}
\end{example}

The concepts of continuity given in Definitions \ref{Defscont},
\ref{Defnormal} and \ref{DefHcont} are strongly connected to the
concepts of semi-continuity of real functions. We have the
following characterization of the fixed points of the lower and
the upper Baire operators, see Appendix 1:
\begin{eqnarray}
I(f)=f&\Longleftrightarrow & f \mbox{ is lower semi-continuous on
}\Omega\label{flsc}\\
S(f)=f & \Longleftrightarrow & f \mbox{ is upper semi-continuous
on }\Omega\label{fusc}%
\end{eqnarray}
Hence for an interval function $f=[\underline{f},\overline{f}]$%
\begin{equation}\label{sclscusc}
f \mbox{ is S-continuous}\;\Longleftrightarrow
\;\left\{\begin{tabular}{c} $\underline{f}$ is lower
semi-continuous\\\\
$\overline{f}$ is upper semi-continuous
\end{tabular}
\right.
\end{equation}
The name of the D-continuous interval functions is due to a
similar characterization through the normal upper semi-continuous
and normal lower semi-continuous functions, see Appendix 1. We
have
\begin{equation}\label{normalnlscnusc}
f \mbox{ is D-continuous}\;\Longleftrightarrow
\;\left\{\begin{tabular}{c} $\underline{f}$ is normal lower
semi-continuous\\\\
$\overline{f}$ is normal upper semi-continuous
\end{tabular}
\right.
\end{equation}
The minimality condition associated with the Hausdorff continuous
functions can also be formulated in terms of semi-continuous
functions, namely, if $f=[\underline{f},\overline{f}]$ is
S-continuous then $f$ is H-continuous if and only if
\begin{equation}\label{mincond}\label{Hcontsemicont}
\{\phi\in {\cal A}(\Omega):\phi\;\mbox{ is
semi-continuous,}\;\underline{f}\leq\phi\leq\overline{f}\}=\{\underline{f},\overline{f}\}
\end{equation}

All three concepts of continuity defined in this section can be
considered as generalizations of the concept of continuity of real
functions in the sense that the only real (point valued) functions
contained in each one of the above sets are the continuous
functions. This is formally stated in the next theorem.

\begin{theorem}\label{thcontcont}
If the function $f\in{\cal A}(\Omega)$ is S-continuous,
D-continuous or H-continuous than it is continuous.
\end{theorem}
\begin{proof}
In view of the implications in Theorem \ref{thcontregscont}, it is
enough to consider the case when $f$ is S-continuous. If
$f=[f,f]\in{\cal A}(\Omega)$ is S-continuous, then $f$ is both
upper semi-continuous and lower semi-continuous on $\Omega$, see
(\ref{sclscusc}). Therefore $f$ is continuous on $\Omega$.
\end{proof}

{\bf Historical Remark.} The three concepts of continuity of
interval functions discussed in this section, namely,
S-continuity, D-continuity and H-continuity, are linked with the
concepts of semi-continuity of real functions, see
(\ref{sclscusc}), (\ref{normalnlscnusc}) and
(\ref{Hcontsemicont}). The lower and upper semi-continuous
functions have been well known at least since the beginning of the
20th century and are usually credited to Baire, see \cite{Baire}.
The normal upper semi-continuous functions were introduced in 1950
by Dilworth in connection with the order completion of the lattice
of continuous functions, see \cite{Dilworth}. The concepts of
S-continuity and H-continuity are both due to Sendov, see
\cite{Sendov0}, \cite{Sendov}. It is quite interesting that
pairing a lower semi-continuous function $\underline{f}$ with an
upper semi-continuous function $\overline{f}$, such that
$\underline{f}\leq\overline{f}$ produces a completely new concept
from both algebraic and topological points of view, namely, the
concept of S-continuous interval functions. It is shown in
\cite{Sendov} that the set of all S-continuous functions on a
compact subset of $\mathbb{R}$ is a complete metric space with
respect to the Hausdorff distance between their graphs and has the
rare and particularly useful property of being completely bounded.
Similarly to the concept of S-continuity, here we consider
interval functions given by pairs of a normal lower
semi-continuous function $\underline{f}$ and a normal upper
semi-continuous function $\overline{f}$, such that
$\underline{f}\leq\overline{f}$ are considered. In honor of the
contribution of Dilworth to this development these functions are
called Dilworth continuous, or shortly D-continuous.

\section{Spaces of Hausdorff continuous interval functions}

The H-continuous functions, representing the strongest kind of
continuity among the three considered above, are also similar to
the usual continuous real functions in that they assume point
(degenerate interval) values on a dense subset of the domain
$\Omega$. This is obtained from a Baire category argument. It was
shown in \cite{QM} that for every $f\in\mathbb{H}(\Omega)$ the set
\begin{equation}\label{Gammaf}
W_f=\{x\in\Omega : w(f(x))>0\}
\end{equation}
is of first Baire category. Since $\Omega\subseteq \mathbb{R}^n$
is open this implies that for every $f\in\mathbb{H}(\Omega)$ the
set
\begin{equation}\label{Dfdense}
D_f=\{x\in\Omega : w(f(x))=0\}=\Omega\backslash W_f\ \;{\rm
is\;dense\;in\;}\Omega.
\end{equation}
Using that a finite or countable union of sets of first Baire
category is also a set of first Baire category we have that for
every finite or countable set ${\cal F}$ of Hausdorff continuous
functions the set
\begin{equation}\label{Dsetfdense}
D_{\cal F}=\{x\in\Omega : w(f(x))=0,\;f\in{\cal
F}\}=\Omega\backslash\left(\displaystyle\bigcup_{f\in{\cal F}}
W_f\right)\;\;{\rm is\;dense\;in\;}\Omega.
\end{equation}

The property (\ref{Dfdense}) can also be used to characterize the
H-continuous functions as stated in the following theorem.

\begin{theorem}\label{trH}
If the interval function $f$ is D-continuous and
assumes point (degenerate interval) values on a dense subset $D$ of $\Omega $, that is,%
\[
w(f(x))=0,\ x\in D,\ D \mbox{ -- dense subset of }\Omega ,
\]%
then $f$ is H-continuous.
\end{theorem}
\begin{proof}Assume that the function $f$ is D-continuous and satisfies
the condition in the theorem.  Let $g\in\mathbb{A}(\Omega)$ be
such that
\[
g(x)\subseteq f(x),\ x\in\Omega.
\]
Then we have
\[
g(x)=f(x),\ x\in D.
\]
Hence
\[
F(g)(x)\subseteq
F(f)(x)=f(x)=F(D,\Omega,f)(x)=F(D,\Omega,g)(x)\subseteq
F(g)(x),x\in \Omega,
\]
where for the last inclusion we use the property (\ref{Finclmon3})
of the operator $F$. The above inclusions indicate that $F(g)=f$
which means that $f$ is H-continuous.
\end{proof}

It may appear at first that the minimality condition in Definition
\ref{DefHcont} applies at each individual point $x$ of $\Omega$,
thus, not involving neighborhoods. However, the graph completion
operator $F$ does appear in this condition. And this operator
according to (\ref{gcgen}) and therefore (\ref{lbfgen}) and
(\ref{ubfgen}) does certainly refer to neighborhoods of points in
$\Omega$, a situation typical, among others, for the concept of
continuity. Hence the following property of the continuous
functions is preserved, \cite{QM}.

\begin{theorem}
\label{tindent}Let $f,g$ be H-continuous on $\Omega$ and let $D$
be a dense subset
of $\Omega$. Then%
\begin{eqnarray*}
&{\rm a)}& \ \ f(x)~\leq~ g(x),\;x\in D~\Longrightarrow~
f(x)~\leq~ g(x),\;x\in \Omega,{\rm  \ \ \ \ \ \ \ \ \ \ \ \ \ \ \
\ \ \ \ \
}\\
&{\rm b)}& \ \ f(x)~=~g(x),\;x\in D~\Longrightarrow~ f(x)~=~ g(x),\;x\in \Omega.{\rm  \ \ \ \ \ \ \ \ \ \ \ \ \ \ \ \ \ \ \ \ }%
\end{eqnarray*}
\end{theorem}

The following two theorems represent essential links with the
usual point valued continuous functions.

\begin{theorem}\label{tcont}
Let $f=[\underline{f},\overline{f}]$ be an H-continuous function
on $\Omega$. \newline a) If $\underline{f}$ or $\overline{f}$ is
continuous at a point $a\in \Omega$ then
$\underline{f}(a)=\overline{f}(a)$.\newline b) If
$\underline{f}(a)=\overline{f}(a)$ for some $a\in \Omega$ then
both $\underline{f}$ and $\overline{f}$ are continuous at $a$.
\end{theorem}
The proof is given in \cite{QM}

\begin{theorem}\label{textend}
Let $D$ be a dense subset of $\Omega$. If $f\in C(D)$ then
\begin{eqnarray*}
&(i)&F(D,\Omega,f)(x)=f(x),\ x\in D,\\
&(ii)&F(D,\Omega,f)\in\mathbb{H}(\Omega),\ \ \ \ \ \ \ \ \ \ \ \ \
\ \ \ \ \ \ \ \ \ \ \ \
 \end{eqnarray*}
\end{theorem}
\begin{proof}
(i) Using that $f$ is continuous on $D$, for every $x\in D$ we
have
\[F(D,\Omega,f)(x)=F(D,D,f)(x)=f(x)\]

(ii) Let $g\in\mathbb{A}(\Omega)$ satisfy the inclusion
\[
g(x)\subseteq F(D,\Omega,f)(x),\ x\in\Omega.
\]
Using the property (\ref{Finclmon1}) of the operator $F$ and its
idempotence, see (\ref{idemp}), we have
\[g(x)\subseteq F(g)(x)\subseteq F(F(D,\Omega,f))(x)=F(D,\Omega,f)(x),x\in \Omega .\]
Therefore
\[g(x)=F(D,\Omega,f)(x)=f(x),x\in D\]
Hence
\[
F(g)(x)\subseteq F(D,\Omega,f)(x)=F(D,\Omega,g)(x)\subseteq
F(g)(x),\ x\in\Omega,
\]
where for the last inclusion we used the property
(\ref{Finclmon3}) of the operator $F$. From the above inclusions
we have
\[
F(g)(x)=F(D,\Omega,f)(x)\ ,
\]
which shows that the function $F(D,\Omega,f)$ is H-continuous.
\end{proof}

We state below one of the most amazing properties of the set
$\mathbb{H}(\Omega)$, namely, its order completeness, as well as
the Dedekind order completeness of some of its important subsets.
What makes this property so significant is the fact that with very
few exceptions the usual spaces in Real Analysis or Functional
Analysis are neither order complete nor Dedekind order complete,
see Appendix 2 for the definitions of order completeness and
Dedekind order completeness.

\begin{theorem}\label{tocomp}
The set $\mathbb{H}(\Omega)$ of all H-continuous interval
functions is order complete, that is, for every subset
$\mathcal{F}$ of $\mathbb{H}(\Omega)$ there exist
$u,v\in\mathbb{H}(\Omega)$ such that $u=\sup\mathcal{F}$ and
$v=\inf \mathcal{F}$.
\end{theorem}
\begin{proof}
First we will construct $u=\sup\mathcal{F}$. Consider the function%
\[
g(x)=\sup\{I(f)(x):f\in\mathcal{F}\},\;x\in\Omega.
\]
The function $g$, being a supremum of lower semi-continuous
functions, is also a lower semi-continuous function, see\ Theorem
\ref{tsupinf} in Appendix 1. Therefore, according to Theorem
\ref{tFSI} we have that $u=F(S(g))=F(S(I(g)))$ is H-continuous,
that is $u\in\mathbb{H}(\Omega)$. We will prove that $u$ is the
supremum of $\mathcal{F}$. More precisely, we will show that

i) $u$ is an upper bound of $\mathcal{F}$, i.e. $f(x)\leq
u(x),\;\;x\in\Omega,\;\;f\in\mathcal{F}$,

ii) $u$ is the smallest upper bound of $\mathcal{F}$, that is, for
any function
$h\in\mathbb{H}(\Omega)$,%
\[
f(x)\leq
h(x),\;\;x\in\Omega,\;\;f\in\mathcal{F}\;\;\Longrightarrow\;\;u(x)\leq
h(x),\;\;x\in\Omega.
\]
Using the monotonicity of the operators $S$ and $F,$ for every
$f\in \mathcal{F}$ we have
\begin{eqnarray*}
f(x) & \leq & S(f)(x)\ =\ S(I(f))(x)\ \leq\  S(g)(x),\;x\in\Omega,\\
f(x) & = & F(f)(x)\ \leq\  F(S(g))(x)\ =\ u(x),\;x\in\Omega,
\end{eqnarray*}
which means that $u$ is an upper bound of the set $\mathcal{F}$.

Assume that function $h\in\mathbb{H}(\Omega)$ is such that
\[
f(x)\leq h(x),\;\;x\in\Omega,\;\;f\in\mathcal{F}.
\]
We have%
\[
g(x)=\sup\{I(f)(x):f\in\mathcal{F}\}\leq h(x),\;x\in\Omega,
\]
which implies $S(g)\leq S(h).$ Hence $u=F(S(g))\leq F(S(h))=h.$
Therefore, $u=\sup\mathcal{F}$.

The existence of $v=\inf\mathcal{F}\in\mathbb{H}(\Omega)$ is shown
in a similar way.
\end{proof}

The following theorem gives a useful representation of the infimum
and supremum of a subset of $\mathbb{H}(\Omega)$ in terms of the
point-wise infimum and supremum, respectively.

\begin{theorem}
\label{tpointsupinf}Let $\mathcal{F}\subseteq\mathbb{H}(\Omega)$
and let the
functions $\varphi,\psi\in\mathcal{A}(\Omega)$ be defined by%
\[
\varphi(x)=\inf\{z\in f(x):f\in\mathcal{F}\}{\rm\ , \ \
}\psi(x)=\sup\{z\in f(x):f\in\mathcal{F}\}\ ,{\rm  \ }x\in\Omega.
\]
Then
\[
\inf\mathcal{F}=F(I(\varphi)){\rm , \ \
}\sup\mathcal{F}=F(S(\psi)).
\]
\end{theorem}

\begin{proof}
We will prove that $\sup\mathcal{F}=F(S(\psi))$. The proof of
$\inf \mathcal{F}=F(I(\varphi))$ can be done in a similar way. In
the proof of Theorem \ref{tocomp} the supremum
$u=\sup\mathcal{F}$\ was constructed in the form $u=F(S(g))$ where
\[
g(x)=\sup\{I(f)(x):f\in\mathcal{F}\},\;x\in\Omega.
\]
From the inequality
\[
I(f)\leq f\leq u,\;f\in\mathcal{F},
\]
it follows that
\[
g\leq\psi\leq u.
\]
Using the monotonicity of the operators $S$ and $F$, see
(\ref{mon}), from
the above inequalities we have%
\[
u=F(S(g))\leq F(S(\psi))\leq F(S(u))=u.
\]
Hence $F(S(\psi))=u=\sup\mathcal{F}$.
\end{proof}

\begin{example}
Let $\Omega=\mathbb{R}^{n}$. Consider the set
\[
\mathcal{F}=\{f_{\delta}:\delta>0\}\subseteq C(\Omega),
\]
where%
\[
f_{\delta}(x)=\left\{
\begin{tabular}
[c]{ll}%
$1-\delta^{-1}||x||$ & if $x\in B_\delta(0)$\\
$0$ & otherwise
\end{tabular}
\ \ \right.
\]
The point-wise infimum of the set $\mathcal{F}$ is
\[
\varphi(x)=\left\{
\begin{tabular}
[c]{ll}%
$1$ & if $x=0$\\
$0$ & if $x\neq0$%
\end{tabular}
\ \ \right.
\]
which is not an H-continuous function. The infimum of
$\mathcal{F}$ in $\mathbb{H}(\mathbb{R}^{n})$ is $u(x)=0$,
$x\in\mathbb{R}^{n}$. Clearly, $u=F(I(\varphi))$.
\end{example}

\begin{example}
Let $\Omega=\mathbb{R}$. Consider the set
$\mathcal{F}=\{f_{n}:n\in\mathbb{N}\}$ where
\[
f_{n}(x)=\left\{
\begin{tabular}
[c]{lll}%
$x^{-2n-1}$ & if & $x\neq0$\\
$\lbrack-\infty,\infty]$ & if & $x=0$%
\end{tabular}
\ \ \ \right.  .
\]
The point-wise supremum of $\mathcal{F}$ is%
\[
\psi(x)=\left\{
\begin{tabular}
[c]{lll}%
0 & if & $x<-1$\\
$x^{-1}$ & if & $-1\leq x<0$\\
$\infty$ & if & $0\leq x<1$\\
$x^{-1}$ & if & $x\geq1$%
\end{tabular}
\ \ \ \right.  .
\]
Hence%
\[
\sup\mathcal{F}=\left\{
\begin{tabular}
[c]{lll}%
0 & if & $x<-1$\\
$\lbrack-1,0]$ & if & $x=-1$\\
$x^{-1}$ & if & $-1<x<0$\\
$\lbrack-\infty,\infty]$ & if & $x=0$\\
$\infty$ & if & $0<x<1$\\
$\lbrack1,\infty]$ & if & $x=1$\\
$x^{-1}$ & if & $x>1$%
\end{tabular}
\ \ \ \right.  .
\]
\end{example}

\begin{remark}
The Theorem \ref{tpointsupinf} establishes a close connection
between the supremum (infimum) in $\mathbb{H}(\Omega)$ with
respect to the partial order (\ref{forder}) and the point-wise
supremum (infimum). However, as the above examples show, these two
functions are clearly not the same, that is, for a set
$\mathcal{F}\subseteq\mathbb{H}(\Omega)$, in general,
\[
\left(  \sup\mathcal{F}\right)  (x) =\left(  \sup_{f\in\mathcal{F}%
}f\right)  (x)\neq\sup_{f\in\mathcal{F}}\left(f(x)\right)\ , \ x\in\Omega{\rm \ .}%
\]
\end{remark}

The next theorem states the Dedekind order completeness of some
important subsets of $\mathbb{H}(\Omega)$.

\begin{theorem}\label{tocomp2}The following subsets of
$\mathbb{H}(\Omega)$are all Dedekind order complete:
\begin{itemize}
\item[(i)] the set $\mathbb{H}_{bd}(\Omega)$ of all bounded
H-continuous functions;

\item[(ii)] The set $\mathbb{H}_{f\!t}(\Omega)$ of all finite
H-continuous functions;

\item[(iii)] The set $\mathbb{H}_{n\!f}(\Omega)$ of all nearly
finite H-continuous.
\end{itemize}
\end{theorem}
The proof of (i) and (ii) can be found in \cite{QM}. Unlike the
concepts of boundedness or finitness associated with the sets
$\mathbb{H}_{bd}(\Omega)$ and $\mathbb{H}_{f\!t}(\Omega)$ which
are well known, the concept of a function being nearly finite is
relatively recent. Further details and properties associated with
this concept are presented in Appendix 3, where the proof of (iii)
can also be found.

\section{The space of Dilworth continuous interval functions}

The D-continuous functions are closely linked to the normalizing
operator $G$ as shown already in Theorem \ref{tFSI}. Hence we will
first state some properties of this operator.

\begin{theorem}\label{tnormop}
The normalizing operator $G$ is

(i)monotone increasing with respect to the partial order
(\ref{forder}), that is, for every two functions
$f,g\in\mathbb{A}(\Omega)$ we have
\[
f(x)\leq g(x),\ x\in\Omega\ \Longrightarrow\ G(f)(x)\leq G(g)(x),\
x\in\Omega;
\]

(ii) inclusion isotone, that is, for every two functions
$f,g\in\mathbb{A}(\Omega)$ we have
\[
f(x)\subseteq g(x),\ x\in\Omega\ \Longrightarrow\ G(f)(x)\subseteq
G(g)(x),\ x\in\Omega;
\]

(iii) idempotent, that is, for every $f\in\mathbb{A}(\Omega)$ we
have
\[
G(G(f))=G(f)
\]
\end{theorem}
\begin{proof}
The statements (i) and (ii) follow directly for the monotonicity
of the operators $I$ and $S$, see (\ref{mon}).

(iii) The normalizing operator can be represented in the form
\[G(f)=[(S\circ I)(f),(I\circ S)(f)],\ f\in\mathbb{A}(\Omega).\]
Therefore, to prove that $G$  is idempotent it is enough to show
that the compositions
\begin{equation}\label{idempIS}
I\circ S\ \ {\rm and}\ \ S\circ I
\end{equation}
are both idempotent. This follows from the idempotence, see
(\ref{idemp}), of the operators $I$ and $S$ and their
monotonicity, see (\ref{mon}). Indeed, for every
$f\in\mathbb{A}(\Omega)$ we have
\begin{eqnarray*}
((I\circ S)\circ(I\circ S))(f)&=&I(S(I(S(f))))~\leq~
I(S(S(S(f))))~=~I(S(f))~=~(I\circ S)(f)\\
((I\circ S)\circ(I\circ S))(f)&=&I(S(I(S(f))))~\geq~
I(I(I(I(f))))~=~I(S(f))~=~(I\circ S)(f)
\end{eqnarray*}
Hence
\[(I\circ S)\circ(I\circ S)=(I\circ S)
\]
The idempotence of $S\circ I$ is proved in the same way.
\end{proof}

From the monotonicity of the operators $I$ and $S$ one can also
easily see that, for every $f\in\mathbb{A}(\Omega)$
\[
G(f)(x)\subseteq F(f)(x),\ x\in\Omega.
\]
Furthermore, we have
\[
G\circ F=F\circ G=G.
\]

The next two theorems give necessary and sufficient conditions for
a function to be D-continuous, the first one - in terms of
H-continuity, the second one - in terms of the operator $G$.

\begin{theorem}\label{tnormalHcont}
A function $f=[\underline{f},\overline{f}]\in\mathbb{A}(\Omega)$
is D-continuous if and only if the functions $F(\underline{f})$
and $F(\overline{f})$ are both H-continuous.
\end{theorem}
\begin{proof} Assume first that the functions
\begin{eqnarray*}
F(\underline{f})&=&[\underline{f},S(\underline{f})]\\
F(\overline{f})&=&[I(\overline{f}),\overline{f}]
\end{eqnarray*}
are H-continuous. From Theorem \ref{thcontregscont} it follows
that these functions are D-continuous as well. Let $D$ be a dense
subset of $\Omega$. Using that $F(\underline{f})$ is D-continuous
we obtain
\[I(D,\Omega,\underline{f})=\underline{f},\]
while using that $F(\overline{f})$ is D-continuous we obtain
\[S(D,\Omega,\overline{f})=\overline{f}
\]
The above two equations put together give
\[
F(D,\Omega,f)=f
\]
which means that $f$ is D-continuous.

Let now $f=[\underline{f},\overline{f}]\in\mathbb{A}(\Omega)$ be
D-continuous. We will show that
$g=[\underline{g},\overline{g}]=F(\underline{f})$ is H-continuous.
Clearly $\underline{g}=\underline{f}$ and
$\overline{g}=S(\underline{f})$. Let $\varepsilon>0$ be fixed.
Since the function
\[
w(g)(x)=w(g(x))=\overline{g}(x)-\underline{g}(x),\ x\in\Omega,
\]
is upper semi-continuous, the set
\[
W_{g,\varepsilon}=\{x\in\Omega:w(g(x))\geq\varepsilon\}
\]
is closed. We will show that this set is nowhere dense as well.
Assume the opposite, that is, there exists an open subset $P$ of
$\Omega$ such that $W_{g,\varepsilon}$ is dense in $P$. Then using
again the fact that the function $w(g)$ is upper semi-continuous,
for every $x\in P$ we have
\begin{eqnarray*}
w(g)(x)&=&S(w(g))(x)=\inf_{\delta>0}\sup\{w(g)(y):y\in
B_\delta(x)\}\\
&\geq&\inf_{\delta>0}\sup\{w(g)(y):y\in B_\delta(x)\cap
W_{g,\varepsilon}\}~\geq~\varepsilon
\end{eqnarray*}
Therefore
\begin{equation}\label{eqwgeps}
S(\underline{f})(x)-\underline{f}(x)=\overline{g}(x)-\underline{g}(x)=w(g)(x)\geq
\varepsilon,\ x\in P.
\end{equation}
Let us fix $x\in P$. Since $P$ is open there exists $\delta_0>0$
such that $B_{\delta_0}(x)\subseteq P$. Now, using also
(\ref{eqwgeps}) we obtain
\begin{eqnarray*}
S(\underline{f})(x)&=&\inf_{\delta>0}\sup\{\underline{f}(y):y\in
B_\delta(x)\}\\
&=&\inf_{0<\delta<\delta_0}\sup\{\underline{f}(y):y\in
B_\delta(x)\}\\
&\leq&\inf_{0<\delta<\delta_0}\sup\{S(\underline{f})(y)-\varepsilon:y\in
B_\delta(x)\}\\
&=&S(S(\underline{f}))-\varepsilon~=~S(\underline{f})-\varepsilon.
\end{eqnarray*}
The obtained contradiction shows that the set $W_{g,\varepsilon}$
is nowhere dense. Then the set
\[
W_g=\bigcup_{\varepsilon>0}W_{g,\varepsilon}
\]
being a union of closed, nowhere dense sets is a set of first
Baire category. This implies that its complement in $\Omega$
\[
D_g=\Omega\setminus
W_g=\{x\in\Omega:S(\underline{f})(x)=\underline{f}(x)\}
\]
is dense in $\Omega$. From the D-continuity of $f$ it follows that
\[
I(D_g,\Omega,\underline{f})=\underline{f}.
\]
Hence
\[
I(S(\underline{f}))\leq
I(D_g,\Omega,S(\underline{f}))=I(D_g,\Omega,\underline{f})=\underline{f}
\]
Since the inequality $I(S(\underline{f}))\geq \underline{f}$ is
obvious we have
\[
I(S(\underline{f}))=\underline{f}
\]
Therefore $g$ can be represented in the form
\[
g=F(\underline{f})=[\underline{f},S(\underline{f})]=[I(S(I(f))),S(I(f))=F(S(I(f)))
\]
and the H-continuity of $g$ follows from Theorem \ref{tFSI}.

The H-continuity of $F(\overline{f})$ is proved in the same way.
\end{proof}

\begin{theorem}
A function $f\in\mathbb{A}(\Omega)$ is D-continuous if and only if
\[
G(f)=f
\]
\end{theorem}
\begin{proof}
Let $f\in\mathbb{A}(\Omega)$ be such that $G(f)=f$. It follows
from Theorem \ref{tFSI} that $G(f)$ is D-continuous. Then $f$ is
D-continuous as well.

Assume now that $f\in\mathbb{A}(\Omega)$ is D-continuous. It
follows from Theorem \ref{tnormalHcont} that the functions
\begin{eqnarray*}
F(\underline{f})&=&[\underline{f},S(\underline{f})]\\
F(\overline{f})&=&[I(\overline{f}),\overline{f}]
\end{eqnarray*}
are H-continuous. The minimality condition (\ref{mincond})
indicates that
\begin{eqnarray*}
I(S(\underline{f}))&=&\underline{f},\\
S(I(\overline{f}))&=&\overline{f}.
\end{eqnarray*}
Thus $G(f)=f$.
\end{proof}

The D-continuous functions are not 'thin' as the H-continuous
functions, that is, they may assume interval values on open
subsets of $\Omega$ and on the whole of $\Omega$ for that matter.
However, they still retain the properties of the continuous and
H-continuous functions stated in Theorem \ref{tindent}, this time
also extended with the relation inclusion as formulated in the
next theorem.

\begin{theorem}
Let $f,g$ be D-continuous on $\Omega$ and let $D$ be a dense
subset
of $\Omega$. Then%
\begin{eqnarray*}
&{\rm a)}& \ \ f(x)~\leq~ g(x),\;x\in D~\Longrightarrow~
f(x)~\leq~ g(x),\;x\in \Omega,{\rm  \ \ \ \ \ \ \ \ \ \ \ \ \ \ \
\ \ \ \ \
}\\
&{\rm b)}& \ \ f(x)~=~g(x),\;x\in D~\Longrightarrow~ f(x)~=~ g(x),\;x\in \Omega.
{\rm  \ \ \ \ \ \ \ \ \ \ \ \ \ \ \ \ \ \ \ \ }\\
&{\rm c)}&\ \ f(x)\subseteq g(x),\ x\in D~\Longrightarrow~
f(x)~\subseteq~ g(x),\;x\in \Omega.
\end{eqnarray*}
\end{theorem}
The proof follows directly from Definition \ref{Defnormal}.

The set of D-continuous functions has similar order completeness
properties to the properties of the set of Hausdorff continuous
functions stated in Theorems \ref{tocomp} and \ref{tocomp2}.

\begin{theorem}\label{tocomplnormal}
(i) The set of $\mathbb{G}(\Omega)$ all H-continuous functions is
order complete with respect to the partial order (\ref{forder})

(ii) Each of the following sets is Dedekind order complete with
respect to the partial order (\ref{forder}):
\begin{itemize}
\item $\mathbb{G}_{bd}(\Omega)$ the set of all bounded
D-continuous interval functions;

\item $\mathbb{G}_{f\!t}(\Omega)$ the set of all finite
D-continuous interval functions;

\item $\mathbb{G}_{n\!f}(\Omega)$ the set of all nearly finite
D-continuous interval functions.
\end{itemize}
\end{theorem}

Note that the use of the relation inclusion (\ref{definclusion})
makes little sense for H-continuous function since, due to the
minimality condition in Definition \ref{DefHcont}, two
H-continuous function are compared with respect to
(\ref{definclusion}) if and only if they are equal. The situation
in the case of D-continuous functions is completely different
since these functions may assume proper interval values on open
subsets of $\Omega$ or on the whole of $\Omega$. It is rather
interesting that all sets considered in Theorem
\ref{tocomplnormal} are also Dedekind order complete with respect
to the relation inclusion given in (\ref{definclusion}).

\begin{theorem}\label{tinclocomp}The sets $\mathbb{G}_{bd}(\Omega)$,
$\mathbb{G}_{f\!t}(\Omega)$, $\mathbb{G}_{n\!f}(\Omega)$ and
$\mathbb{G}(\Omega)$ are all Dedekind order complete with respect
to order relation inclusion given in (\ref{definclusion}).
\end{theorem}

The proofs of both Theorem \ref{tocomplnormal} and Theorem
\ref{tinclocomp} follow from Theorem \ref{tocomp} by using that
every D-continuous function can be represented through two
H-continuous functions as provided by Theorem \ref{tnormalHcont}.

\section{The space of S-continuous interval functions}

The space $\mathbb{F}(\Omega)$ is much wider then the set of
D-continuous functions. However, as demonstrated in the preceding
sections, it is a useful embedding structure for the discussed
spaces. In particular, one may note that it contains the completed
graphs of all point-wise infima and suprema of sets of continuous
functions. The space of all finite S-continuous function, which we
denote by $\mathbb{F}_{f\!t}(\Omega)$ is studied in \cite{Sendov}
in case when $\Omega$ is a compact real interval. It was shown
that this space has two very interesting properties when
considered as a metric space with respect to the Hausdorff
distance, namely,
\begin{itemize}
\item $\mathbb{F}_{f\!t}(\Omega)$ is a complete metric
space%
\item the closed and bounded subsets of
$\mathbb{F}_{f\!t}(\Omega)$ are compact.
\end{itemize}
These properties can also be obtained as a consequence of a more
general theorem about the Hausdorff metric on the set of compact
subsets of $\mathbb{R}^n$, see \cite{Brunn-Minkowski}, Theorems
1.8.2 and 1.8.3. The second property is particularly interesting
in view of the fact that it cannot be attributed to the usual
spaces of functions considered in Real Analysis or Functional
Analysis. Similar properties can be formulated for the case when
$\Omega$ is an open subset of $\mathbb{R}^n$, this however is
beyond the scope of this paper and will be discussed  in a
separate publication.

\section{Application: Dedekind order completion of sets of continuous
functions}%
Here we consider the set $C(\Omega)$\ of all continuous
real functions defined on $\Omega$, that is,
\[
C(\Omega)=\{f:X\rightarrow\mathbb{R},\;f{\rm -continuous\ on\
}\Omega\}
\]
with a partial order defined point-wise as in (\ref{forder}), this
time the values of the functions being real numbers. Clearly
$C(\Omega)$ is a subset of $\mathbb{H}(\Omega)$, see Theorem
\ref{thcontcont}. Furthermore, since the order
between intervals (\ref{iorder}) is an extension of the order\ in $\mathbb{R}%
$, the order in $\mathbb{H}(\Omega)$ is an extension of the order
in $C(\Omega)$. In this way $C(\Omega)$ is embedded in the order
complete set $\mathbb{H}(\Omega)$. We will also consider the
subset $C_{bd}(\Omega)$ of $C(\Omega)$ consisting of all bounded
continuous functions, that is,%
\[
C_{bd}(\Omega)=\{f\in C(\Omega):\exists M\in R:|f(x)|\leq M,\;x\in
\Omega\}.
\]
The fact that both $C_{bd}(\Omega)$ and $C(\Omega)$ are neither
order complete nor Dedekind order complete is well known and can
be shown by trivial examples, see Appendix 2 for the respective
definitions. A general result on Dedekind order completion of
partially ordered sets was established by MacNeilly in 1937, see
\cite{Luxemburg} for a more recent presentation. The problem of
order completion of $C(\Omega)$ is particularly addressed in
\cite{Mack}. In our approach we seek characterization of the
Dedekind order completions of $C_{bd}(\Omega)$ and $C(\Omega)$ as
subsets of $\mathbb{H}(\Omega)$. Obviously, since the set
$\mathbb{H}(\Omega)$ is order complete it also contains the
Dedekind order completions of these sets. An interesting aspect of
this approach is that the Dedekind order completions of the
mentioned sets of continuous functions are constructed again as
sets of functions defined on the same domain $\Omega$. An earlier
result of Dilworth is of similar nature. In \cite{Dilworth} it is
proved that the set of the normal upper semi-continuous real
functions defined on $\Omega$ is a Dedekind order completion of
$C_{bd}(\Omega)$. In \cite{QM} an \emph{alternative}
characterization of the Dedekind order completion of
$C_{bd}(\Omega)$ was found and obtained again as a set of
functions, this time given as a subset of the set
$\mathbb{H}(\Omega)$. Furthermore, the Dedekind order completion
of $C(\Omega)$ is also characterized as a subset of
$\mathbb{H}(\Omega)$. More precisely, in terms of the notations
adopted in this paper, we have
\begin{itemize}
\item $\mathbb{H}_{bd}(\Omega)$ is a Dedekind order completion of
$C_{bd}(\Omega)$;

\item $\mathbb{H}_{f\!t}(\Omega)$ is a Dedekind order completion
of $C(\Omega)$.
\end{itemize}

Following the method used to prove the second statement above, see
Theorem 10 in \cite{QM} one can show that the following
representation of H-continuous functions.

\begin{theorem}\label{tocompo} Let $f=[\underline{f},\overline{f}]\in\mathbb{H}(\Omega)$.
\begin{itemize}\item[(i)] If $f(x)>-\infty$, $x\in\Omega$ then
\[
f=\sup\{g\in C(\Omega):g\leq f\}
\]
\item[(ii)]If $f(x)<+\infty$, $x\in\Omega$ then
\[
f=\inf\{g\in C(\Omega):g\geq f\}
\]
\end{itemize}
\end{theorem}

We should note that the representation of the H-continuous
functions given in the above theorem essentially uses Dedekind
cuts. Under the conditions considered in this theorem the
point-wise representation given in Theorem \ref{tpointsupinf}
simplifies in the following way:
\begin{eqnarray}
{\rm case}\ (i):&&\underline{f}(x)=\sup\{g(x):g\in
C(\Omega):g\leq f\},\ x\in\Omega;\label{tocomporep1}\\
{\rm case}\ (ii):&&\overline{f}(x)=\inf\{g(x):g\in C(\Omega):g\leq
f\},\ x\in\Omega;\label{tocomporep2}
\end{eqnarray}

Obviously, the Theorem \ref{tocompo} does not give representation
in terms of functions in $C(\Omega)$ for all H-continuous
functions since there are H-continuous functions with values
involving both $-\infty$ and $+\infty$. However, using Theorem
\ref{tocompo} it is not difficult to prove:

\begin{theorem}
For every $f\in\mathbb{H}(\Omega)$ we have
\[
f=\inf\{\sup\{g\in C(\Omega):g\leq \sup\{f,m\}\}:m\in\mathbb{R}\}
\]
\end{theorem}
In the above theorem the letter $m$ denotes both the real number
$m$ and the constant function with value $m$ on $\Omega$. This
theorem shows that $\mathbb{H}(\Omega)$ is the minimal order
complete set containing $C(X)$, that is,
\begin{itemize}
\item[] $\mathbb{H}(\Omega)$ is an order completion of
$C(\Omega)$.
\end{itemize}

\section{Application: Order isomorphic representation of
piece-wise smooth functions}%
In this section we consider the space $C_{nd}(\Omega)$ of real
(point) valued functions which are continuous every where on
$\Omega$ except for a closed, nowhere dense subset
$\Gamma\subseteq\Omega$. Obviously every possible kind of
piece-wise continuous functions are included in the set
$C_{nd}(\Omega)$ because the set $\Gamma$ may have arbitrary
shapes.  In addition $\Gamma$ may also have a positive Lebesgue
measure. The space $C_{nd}(\Omega)$ is used, among others, with
the order completion method for solution of nonlinear PDEs, see
\cite{Rosinger}.

More precisely the set $C_{nd}(\Omega)$ is defined as
\begin{equation}\label{cndomega}
C_{nd}(\Omega)=\left\{  f\;\left|
\begin{tabular}
[c]{l}%
$\exists\Gamma\subset\Omega$ closed, nowhere dense:\\
i) $f:\Omega\backslash\Gamma\longmapsto\mathbb{R}$\\
ii) $f\in C(\Omega\backslash\Gamma)$%
\end{tabular}
\ \ \ \ \ \ \ \right.  \right\}   %
\end{equation}
Since the only point valued functions in the set
$\mathbb{H}(\Omega)$ are the functions which are continuous on the
whole of $\Omega$, see Theorem \ref{thcontcont}, the set
$C_{nd}(\Omega)$ is not a subset of $\mathbb{H}(\Omega)$. We shall
define a mapping from $C_{nd}(\Omega)$ to $\mathbb{H}(\Omega)$
which gives are representation of the functions in
$C_{nd}(\Omega)$ through Hausdorff continuous functions. To this
end we proceed as follows.

Let $u\in C_{nd}(\Omega)$. According to (\ref{cndomega}), there
exists a closed, nowhere dense set $\Gamma\subset\Omega$ such that
$u\in C(\Omega \setminus\Gamma)$. Since $\Omega\setminus\Gamma$
is open and dense in $\Omega,$ we can define%
\begin{equation}\label{F0def}
F_{0}(u)=F(\Omega\setminus\Gamma,\Omega,u)
\end{equation}
The closed, nowhere dense set $\Gamma$ used in (\ref{F0def}), is
not unique. However, we can show that the value of
$F(\Omega\setminus\Gamma,\Omega,u)$ does not depend on the set
$\Gamma$ in the sense that for every closed, nowhere dense set
$\Gamma$ such that $u\in C(\Omega\setminus\Gamma)$ the value of
$F(\Omega\setminus\Gamma,\Omega,u)$ remains the same.

Let $\Gamma_{1}$ and $\Gamma_{2}$ be closed, nowhere dense sets
such that $u\in C(\Omega\setminus\Gamma_{1})$ and $u\in
C(\Omega\setminus\Gamma_{2})$. Then the set $\Gamma_1\cup\Gamma_2$
is also closed and nowhere dense. According to Theorem
\ref{textend} the functions $F(\Omega\setminus\Gamma_1,\Omega,u)$,
$F(\Omega\setminus\Gamma_2,\Omega,u)$ and
$F(\Omega\setminus(\Gamma_1\cup\Gamma_2),\Omega,u)$ are all
H-continuous and for every $x\in\Omega\setminus(\Gamma_1\cup\Gamma_2)$ we have%
\[
F(\Omega\setminus\Gamma_1,\Omega,u)(x)=
F(\Omega\setminus\Gamma_2,\Omega,u)(x)=
F(\Omega\setminus(\Gamma_1\cup\Gamma_2),\Omega,u)(x)=u(x).
\]
Since $\Omega\setminus(\Gamma_1\cup\Gamma_2)$ is dense in
$\Omega$, Theorem \ref{tindent} implies that %
\[
F(\Omega\setminus\Gamma_1,\Omega,u)=
F(\Omega\setminus\Gamma_2,\Omega,u)=
F(\Omega\setminus(\Gamma_1\cup\Gamma_2),\Omega,u).
\]
Therefore, the mapping
\[F_{0}:C_{nd}(\Omega)\longmapsto\mathbb{A}(\Omega)\]
is unambiguously defined through (\ref{F0def}). In analogy with
(\ref{gcgen}), we call $F_0$ a graph completion mapping on
$C_{nd}(\Omega)$. As mentioned already above it follows from
Theorem \ref{textend} that for every $u\in C_{nd}(\Omega)$ we have%
\[F_0(u)\in\mathbb{H}(\Omega).
\]
Furthermore, if $u\in C(\Omega\setminus\Gamma)$ we have%
\begin{equation}\label{F0inden}
F_0(u)(x)=u(x),\ x\in\Omega\setminus \Gamma.
\end{equation}
The above identity shows that the values of the function $F_0(u)$
are finite on the open and dense set $\Omega\setminus \Gamma$.
Hence, $F_0(u)$
is nearly finite, see Definition \ref{Defnearlyfinite}. Thus, we have%
\begin{equation}\label{CndHnf}
F_{0}:C_{nd}(\Omega)\longmapsto\mathbb{H}_{nf}(\Omega).%
\end{equation}

The following theorem shows that a function $f\in C_{nd}(\Omega)$
can be identified with $F_0(f)$ in a very direct way.

\begin{theorem} Let $f\in C_{nd}(\Omega)$.
\begin{itemize}
\item[(i)] If $f\in C(\Omega\setminus\Gamma)$ then
\[
f(x)=F_0(f)(x),\ x\in\Omega\setminus\Gamma;
\]

\item[(ii)] If $W_{F_0(f)}$ is the subset of $\Omega$ defined
through (\ref{Gammaf}) for the function $F_0(f)$ then
\[
f(x)=F_0(f)(x),\ x\in\Omega\setminus W_{F_0(f)}.
\]
\end{itemize}
\end{theorem}

The above theorem shows that
\begin{itemize}
\item the largest set on which $f\in C_{nd}(\Omega)$ can be
defined in a continuous way is $\Omega\setminus W_{F_0(f)}$, that
is, for every set $\Gamma$ which is associated with $f$ in terms
of (\ref{cndomega}) we have $W_{F_0(f)}\subseteq\Gamma$;

\item the set $W_{F_0(f)}$ is a closed, nowhere dense set and not
merely a set of first Baire category.

\end{itemize}

Therefore, every function $f\in C_{nd}(\Omega)$ can be identified
with $F_0(f)$ in the following way:

\begin{itemize}

\item If the set $\Gamma\subseteq\Omega$ is associated with $f$ in
terms of (\ref{cndomega}) then $f$ is defined and continuous on the set%
\[
\Omega\setminus\Gamma\subseteq\Omega\setminus W_{F_0(f)}
\]
\item $f$ can be produced in a continuous way on the open and
dense set $\Omega\setminus W_{F_0(f)}$ and it is identical with
${F_0(f)}$ on this set;

\item $f$ can not be produced in a continuous way of a subset of
$\Omega$ which is larger than $W_{F_0(f)}$.

\end{itemize}
Further research  is aimed at applying the results discussed in
this section to improve the regularity results associated with the
order completion method in \cite{Rosinger}.

\section{Application: The set of D-continuous functions contains
all interval hulls of subsets of $C(\Omega)$}

It was shown in Example \ref{ExHullNotCont} that the interval hull
of a set of continuous functions is not always a continuous
function. On the other side, from the representation (\ref{hull})
it is easy to see that the graph of the hull is a closed subset of
$\Omega\times\overline{\mathbb{R}}$ implying that the interval
hull is an S-continuous function. The following theorem gives more
precise characterization.

\begin{theorem}\label{thullnormal}
For every set of continuous real functions $\mathcal{F}\subseteq
C(\Omega)$ the interval hull defined through (\ref{hull}) is a
D-continuous interval function.
\end{theorem}
\begin{proof}
For simplicity we will present the proof only in the case when the
set $\mathcal{F}$ has a finite continuous enclosure.

Let $g=[\underline{g},\overline{g}]={\rm hull}(\mathcal{F})$.
Denote
\[
\mathcal{F}_1=\{\varphi\in C(\Omega):\varphi\leq f,\ \forall
f\in\mathcal{F}\}
\]
Consider $\mathcal{F}$ as a subset of the order complete set
$\mathbb{H}(\Omega)$ and let $\inf\mathcal{F}$ and
$\sup\mathcal{F}$ be the infimum and the supremum of
$\mathcal{F}$, respectively. Then
\begin{equation}\label{F1varphi}
\mathcal{F}_1=\{\varphi\in C(\Omega):\varphi\leq\inf\mathcal{F}\}
\end{equation}
From Theorem \ref{tocompo}(i) we obtain
\begin{equation}\label{infFsupF1}
\inf\mathcal{F}=\sup\mathcal{F}_1
\end{equation}
It follows from (\ref{hull}) that $\underline{g}$ can be
represented as
\begin{equation}\label{gxsup}
\underline{g}(x)=\sup\{\varphi(x):\varphi\in \mathcal{F}_1\},\
x\in\Omega.
\end{equation}
Then using Theorem \ref{tpointsupinf} we have
\[
\sup\mathcal{F}_1=F(S(\underline{g}))
\]
Hence
\[
\inf\mathcal{F}=F(S(\underline{g}))=[I(S(\underline{g})),S(\underline{g})]
\]
Now using (\ref{tocomporep1}) we obtain
\begin{eqnarray*}
I(S(\underline{g}))(x)&=&\sup\{\varphi(x):\varphi\in
C(\Omega):\varphi\leq\inf\mathcal{F}\}\\
&=&\sup\{\varphi(x):\varphi\in \mathcal{F}_1\}\\
&=&\underline{g}(x),\ \ x\in\Omega
\end{eqnarray*}
Therefore
\[
F(\underline{g})=F(S(\underline{g}))=\inf\mathcal{F}\in\mathbb{H}(\Omega).
\]
In the same way we prove that
$F(\overline{g})\in\mathbb{H}(\Omega)$. Then Theorem
\ref{tnormalHcont} implies that $g$ is D-continuous.
\end{proof}

Theorem \ref{thullnormal} shows that the set of D-continuous
functions contains all the interval hulls of all sets of usual
real (point) valued continuous functions. However, in general, the
set $\mathbb{G}_{hl}(\Omega)$ of the interval hulls of sets of
continuous functions is only a subset of $\mathbb{G}(\Omega)$.
More precisely, we have
\[
\mathbb{G}_{hl}(\Omega)=\{f\in\mathbb{G}(\Omega):\exists
\varphi:\Omega\mapsto\overline{\mathbb{R}}:\varphi \ {\rm is\
continuous\ and}\ \varphi\subseteq f\}
\]
Hence, $\mathbb{G}_{hl}(\Omega)$ essentially excludes the
H-continuous functions. Indeed, the only H-continuous functions in
the set $\mathbb{G}_{hl}(\Omega)$ are the continuous point valued
functions on $\Omega$. More precise analysis reveals that the set
of all D-continuous functions functions actually consists of the
interval hulls of the sets of H-continuous functions.

\section*{Appendix 1: Semi-continuous functions and Baire operators}

We recall here the definitions of lower and upper semi-continuity
as given in \cite{Baire}, which also include functions with
extended real values.

\begin{definition}
\label{dlsc}A function $f\in\mathcal{A}(\Omega)$ is called lower
semi-continuous at $x\in\Omega$ if for every $m<f(x)$ there exists
$\delta>0$ such that $m<f(y)$ for all $y\in B_\delta (x).$ If
$f(x)=-\infty$, then $f$ is assumed lower semi-continuous at $x$.
\end{definition}

\begin{definition}
\label{dusc}A function $f\in\mathcal{A}(\Omega)$ is called upper
semi-continuous at $x\in\Omega$ if for every $m>f(x)$ there exists
$\delta>0$ such that $m>f(y)$ for all $y\in B_\delta (x).$ If
$f(x)=+\infty$, then $f$ is assumed upper semi-continuous at $x$.
\end{definition}

\begin{definition}
\label{dluscg}A function $f\in\mathcal{A}(\Omega)$ is called lower
(upper) semi-continuous on $\Omega$ if it is lower (upper)
semi-continuous at every point of $\Omega$.
\end{definition}

The next theorem was proved in \cite{Baire}.

\begin{theorem}
\label{tsupinf} We have the following:\newline a) Let $L\subseteq
\mathcal{A}(\Omega)$ be a set of lower semi-continuous functions.
Then function $l$ defined by
\[
l(x)=\sup\{f(x):f\in L\}
\]
is lower semi-continuous.\smallskip\newline b) Let
$U\subseteq\mathcal{A}(\Omega)$ be a set of upper semi-continuous
functions. Then function $u$ defined by
\[
u(x)=\inf\{f(x):f\in U\}
\]
is upper semi-continuous.
\end{theorem}

From the definitions of the lower and upper Baire operators given
in (\ref{lbfgen}) and (\ref{ubfgen}) and the above theorem one can
immediately see that for every dense subset $D$ of $\Omega$ and
$f\in\mathbb{A}(D)$ we have
\begin{itemize}
\item$I(D,\Omega,f)$ is lower semi-continuous on $\Omega$; \item $
S(D,\Omega,f)$ is upper semi-continuous on $\Omega$.
\end{itemize}
Furthermore, for a given dense subset $D$ of $\Omega$ and
$f\in\mathbb{A}(D)$  the functions $I(D,\Omega,f)$ and
$S(D,\Omega,f)$ have the following optimality properties. For any
$g\in\mathbb{A}(\Omega)$
\begin{eqnarray*}
\left.\begin{tabular}{c}
$g$- lower semi-continuous on $\Omega$\\
$g(x)\leq f(x),\;x\in D$
\end{tabular}
\right\}&\Longrightarrow& g(x)\leq
I(D,\Omega,f)(x),\;x\in\Omega;\\\\
\left.\begin{tabular}{c}
$g$- upper semi-continuous on $\Omega$\\
$g(x)\geq f(x),\;x\in D$
\end{tabular}
\right\}&\Longrightarrow& g(x)\geq S(D,\Omega,f)(x),\;x\in\Omega.
\end{eqnarray*}
Due to the above properties the functions $I(D,\Omega,f)$ and
$S(D,\Omega,f)$ are also called respectively lower and upper
semi-continuous envelops of the function $f$, see \cite{Bardi}.

The following two concepts were introduced by Dilworth,
\cite{Dilworth}.

\begin{definition}
\label{dnlsc}A function $f\in\mathcal{A}(\Omega)$ is called normal
lower semi-continuous on $\Omega$ if it is lower semi-continuous
and
\[I(S(f))=f\ .\]
\end{definition}

\begin{definition}
\label{dnusc}A function $f\in\mathcal{A}(\Omega)$ is called upper
semi-continuous at if it is upper semi-continuous and
\[S(I(f))=f\ .\]
\end{definition}

An important property of the normal lower and upper
semi-continuous is that they can be represented through Dedekind
cuts of the set of continuous functions. More precisely for a
normal lower semi-continuous $f$ we have
\[
f(x)=\sup\{\varphi(x):\varphi\in\mathcal{F}\},\ x\in\Omega,
\]
where
\begin{eqnarray*}
\mathcal{F}&=&\{\varphi:\Omega\mapsto\overline{\mathbb{R}}:\varphi\
{\rm is\ continuous\ on}\ \Omega\ {\rm and}\ \varphi\leq \phi,\ \forall \phi\in\mathcal{G}\}\\
\mathcal{G}&=&\{\phi:\Omega\mapsto\overline{\mathbb{R}}:\phi\ {\rm
is\ continuous\ on}\ \Omega\ {\rm and}\ \phi\geq f\}.
\end{eqnarray*}
In a similar way if $f$ is normal upper semi-continuous then
\[
f(x)=\inf\{\varphi(x):\varphi\in\mathcal{F}\},\ x\in\Omega,
\]
where
\begin{eqnarray*}
\mathcal{F}&=&\{\varphi:\Omega\mapsto\overline{\mathbb{R}}:\varphi\
{\rm is\ continuous\ on}\ \Omega\ {\rm and}\ \varphi\geq \phi,\ \forall \phi\in\mathcal{G}\}\\
\mathcal{G}&=&\{\phi:\Omega\mapsto\overline{\mathbb{R}}:\phi\ {\rm
is\ continuous\ on}\ \Omega\ {\rm and}\ \phi\leq f\}.
\end{eqnarray*}

\section*{Appendix 2: Partial orders for intervals and interval functions}

Several partial orders have historically been associated with the
set $\mathbb{I\,}\overline{\mathbb{R}}$, namely, \newline $(i)$ \
the inclusion
relation $[\underline{a},\overline{a}]\subseteq\lbrack\underline{b}%
,\overline{b}]\Longleftrightarrow\underline{b}\leq\underline{a}\leq
\overline{a}\leq\overline{b}$\newline $(ii)$ \ the ''strong''
partial order
$[\underline{a},\overline{a}]\preceq\lbrack\underline{b},\overline
{b}]\Longleftrightarrow\overline{a}\leq\underline{b}$\newline
$(iii)$ \ the partial order defined by (\ref{iorder}).\newline The
use of the inclusion relation on the set
$\mathbb{I\,}\overline{\mathbb{R}}$ is motivated by the
applications of interval analysis to generating enclosures of
solution sets. However, the role of partial orders extending the
total order on $\overline{\mathbb{R}}$ has also been recognized in
computing, see \cite{Birkhoff}. Both orders $(ii)$ and $(iii)$ are
extensions of the order on $\overline{\mathbb{R}}$. The use of the
order $(ii)$ is based on the view point that inequality between
intervals should imply inequality between their interiors. This
approach is rather limiting since the order $(ii)$ does not retain
some essential properties of the order on $\overline{\mathbb{R}}$.
For instance, a nondegenerate interval $A$ and the interval
$A+\varepsilon$ are not comparable with respect to the order
$(ii)$ when the positive real number $\varepsilon$ is small
enough. The partial order $(iii)$ is introduced and studied by
Markov, see \cite{Markov}, \cite{Markov2}. The results reported in
this paper indicate that indeed the partial order (\ref{forder})
induced point-wise by (\ref{iorder}) is an appropriate partial
order to be associated with the considered spaces of interval
functions. The partial order (\ref{definclusion}) induced
point-wise by the inclusion relation (i) also plays important role
in the spaces where the functions are not 'thin' in the sense that
they assume proper interval values on open subsets of the domain,
e.g. the D-continuous or the S-continuous functions. The following
concepts of order completeness and Dedekind order completeness
were discussed in connection with both orders, namely
(\ref{forder}) and (\ref{definclusion}).

\begin{definition}
\label{defordercomp}A partially ordered set $P$ is called order
complete if every nonempty subset $A$ of $P$ has both a supremum
in $P$ and an infimum in $P.$
\end{definition}

\begin{definition}
\label{defdedordercomp}A partially ordered set $P$ is called
Dedekind order complete if every nonempty subset $A$ of $P$ which
is bounded from above has a supremum in $P$ and every nonempty
subset $B$ of $P$ which is bounded from below has an infimum in
$P.$
\end{definition}

\begin{definition}
\label{defocompletion}Let $P$ be a  partially ordered set. A
partially ordered set $P^{\#}$ is called a (Dedekind) order
completion of $P$ if\newline i) $P^{\#}$ is (Dedekind) order
complete;\newline ii) there exists an order isomorphism
$\Phi:P\rightarrow\Phi(P)\subseteq P^{\#};$\newline iii) if $Q$ is
(Dedekind) order complete and $\Phi(P)\subseteq Q\subseteq P^{\#}$
then $Q=P^{\#}$.
\end{definition}

Clearly, a partially ordered set may be order complete only if all
its subsets are bounded. Let us note here that all subsets of
$\mathbb{H}(\Omega)$ are bounded with respect to the relation
$\leq$. For example $\upsilon(x)=-\infty$, $x\in\Omega$, and
$u(x)=+\infty$, $x\in\Omega$, are, respectively, lower and upper
bounds of every subset of $\mathbb{H}(\Omega)$. However, this is
not a property which is necessarily inherited by the subsets of
$\mathbb{H}(\Omega)$, more precisely, if
$\mathcal{B}\subseteq\mathbb{H}(\Omega)$
then the subsets of $\mathcal{B}$ are not necessarily bounded in $\mathcal{B}%
$.

One should note that the set of Hausdorff continuous function is
Dedekind order complete with respect to the order relation
inclusion ($\subseteq$). However, this statement hardly contains
any information as the only bounded sets with respect to
inclusions are the single function sets. Trivially each of these
sets has both infimum and supremum equal to the function in the
set. Hence we do not have any interest in inclusion within the
space of the Hausdorff continuous functions. As shown in Section 7
the inclusion is an interesting partial order in the wider space
of D-continuous functions.

\section*{Appendix 3: The set of nearly finite H-continuous functions}
With every function $f\in\mathbb{A}(\Omega)$ we associate the set%
\begin{equation}\label{Gammanf}
\Gamma_{n\!f}(f)=\{x\in\Omega:+\infty\in f(x)\;\;{\rm
or}\;\;-\infty\in f(x)\}.
\end{equation}
Then,%
\begin{equation}
f\mbox{ is nearly finite
}\Longleftrightarrow\Gamma_{n\!f}(f)\mbox{ is closed and nowhere dense in }\Omega. \label{iden25}%
\end{equation}

The condition in Definition \ref{Defnearlyfinite}
\emph{simplifies} in the case of H-continuous function as follows.

\begin{theorem}
\label{tcrit} For an H-continuous function $f$ to be nearly finite
it is sufficient to have finite values on a dense subset of
$\Omega$\ which need not be open as well.
\end{theorem}

\begin{proof}
Let us assume that the function $f$ assumes finite values on a set
$D$, which is a dense subset of $\Omega$. According to
(\ref{iden25}) we need to prove that $\Gamma_{n\!f}(f)$ is closed
and nowhere dense in $\Omega$.

Assume that $\Gamma_{n\!f}(f)$ is not nowhere dense. Let
$\Gamma_{+\infty }(f)=\{x\in\Omega:+\infty\in f(x)\}$ and
$\Gamma_{-\infty}(f)=\{x\in\Omega:-\infty\in
f(x)\}$. Clearly $\Gamma_{n\!f}(f)=\Gamma_{+\infty}(f)\cup\Gamma_{-\infty}%
(f)$. Therefore at least one of the sets $\Gamma_{-\infty}(f)$ or
$\Gamma_{+\infty}(f)$ is not nowhere dense, because the union of
two nowhere dense sets is also nowhere dense. Let
$\Gamma_{+\infty}(f)$ be not nowhere dense. Then, there exists an
open set $G\subseteq X$ such that $\Gamma _{+\infty}(f)\cap G$ is
dense in $G.$ The function $f$ assumes finite values on the set
$D$ which is dense in $\Omega$. Hence $G\cap D\neq\emptyset.$ Let
$x\in G\cap D$. Using that $\Gamma_{+\infty}(f)\cap G$ is dense in
$G$ we obtain that for every $\delta>0$ the intersection
$B_\delta(x)\cap\Gamma_{+\infty}(f)$ is not empty. This implies%
\[
\sup\{z\in f(y):y\in B_\delta(x)\}=+\infty
\]
for every $\delta>0$. Therefore%
\[
S(f)(x)=\inf_{\delta>0}\sup\{z\in f(y):y\in B_\delta(x)\}=+\infty.
\]
On the other hand, since $x\in D$ we have $|f(x)|<+\infty$. The
obtained contradiction shows that the assumption is false, i.e.
$\Gamma_{n\!f}(f)$ is nowhere dense in $\Omega$.

We will prove that $\Gamma_{n\!f}(f)$ is closed by proving that
$\Omega\backslash \Gamma_{n\!f}(f)$ is open. Let
$a\in\Omega\backslash\Gamma_{n\!f}(f)$. Then
$I(f)(a)\in\mathbb{R}$ and $S(f)(a)\in\mathbb{R}$. Since $I(f)$
and $S(f)$ are lower and upper semi-continuous functions,
respectively, there exists
$\delta>0$ such that for every $x\in B_\delta(a)$ we have%
\[
I(f)(a)-1\leq I(f)(x)\leq f(x)\leq S(f)(x)\leq S(f)(a)+1.
\]
Therefore $f(x)$ is finite for all $x\in B_\delta(a)$. Hence%
\[
a\in B_\delta(a)\subseteq \Omega\backslash\Gamma_{n\!f}(f).
\]
Thus, $\Omega\backslash\Gamma_{n\!f}(f)$ is an open set.
 \end{proof}

The result in Theorem \ref{tcrit} is further detailed in its
consequences as follows.

\begin{theorem}
\label{thnf}An H-continuous function $f$ which is nearly finite
has the additional property that its values are finite real
numbers, that is, finite point intervals, everywhere on $\Omega$,
except for a set of first Baire category.
\end{theorem}

\begin{proof}
An H-continuous function $f$ assumes proper interval values only
on the set $W_f$, given in (\ref{Gammaf}), which is of first Baire
category, i.e. it is a countable union of closed, nowhere dense
sets. On the other side, the function $f$ assumes nonfinite values
only on the set $\Gamma_{n\!f}(f)$ defined through
(\ref{Gammanf}), which is closed and nowhere dense, see
(\ref{iden25}). Therefore, the function $f$ assumes finite real
values on the set $\Omega\backslash(W_f\cup\Gamma_{n\!f}(f))$. The
set $W_f\cup \Gamma_{n\!f}(f)$ is of first Baire category, because
it is a union of countably many closed, nowhere dense sets. This
completes the proof.
\end{proof}

We consider the set $\mathbb{H}_{n\!f}(\Omega)$ of all
H-continuous nearly finite functions defined on $\Omega$. The
following theorems show that the set $\mathbb{H}_{n\!f}(\Omega)$
is Dedekind order complete with respect to the partial order
defined by (\ref{forder}). Here we should note that a subset
$\mathcal{F}$ of $\mathbb{H}_{n\!f}(\Omega)$ which is bounded from
above or below may still contain functions with values $+\infty$
or $-\infty$ or unbounded closed intervals, this being compatible
with the partial order defined by (\ref{forder}).

\begin{theorem}
\label{tocomp1}The set $\mathbb{H}_{n\!f}(\Omega)$ is Dedekind
order complete, that is,\newline a) if $\mathcal{F}$ is a nonempty
subset of $\mathbb{H}_{n\!f}(\Omega)$ which is bounded from above
then there exists $u\in\mathbb{H}_{n\!f}(\Omega)$ such that
$u=\sup\mathcal{F}$;\newline b) if $\mathcal{F}$ is a nonempty
subset of $\mathbb{H}_{n\!f}(\Omega)$ which is bounded from below
then there exists $v\in\mathbb{H}_{n\!f}(\Omega)$ such that
$v=\inf\mathcal{F}$.
\end{theorem}

\begin{proof}
We will prove only point a). Point b) is proved in a similar way.
Since $\mathcal{F}$ is also a subset of the order complete space
$\mathbb{H}(\Omega)$ the function
$u=\sup\mathcal{F}\in\mathbb{H}(\Omega)$\ is well defined. We only
need to show that $u\in\mathbb{H}_{n\!f}(\Omega)$. We have already
that $u$ is H-continuous. Now we will show that it is nearly
finite. Let $p\in \mathbb{H}_{n\!f}(\Omega)$ be an upper bound of
$\mathcal{F}$. Clearly $u\leq p.$
Let $f\in\mathcal{F}$. Since $u$ is an upper bound of $\mathcal{F}\ $we have%
\[
f(x)\leq u(x)\leq p(x),\;x\in\Omega.
\]
Therefore $u$ assumes finite values at all points of the set
$\Omega\backslash (\Gamma_{n\!f}(f)\cup\Gamma_{n\!f}(u))$, which
is open and dense in $\Omega$. Hence
$u\in\mathbb{H}_{n\!f}(\Omega)$.
\end{proof}


\begin{thebibliography}{00}


\bibitem{QM} R Anguelov, Dedekind Order Completion of C(X) by
Hausdorff Continuous Functions, Quaestiones Mathematicae, to
appear.

\bibitem {Anguelov-Markov}R. Anguelov and S. Markov, Extended segment
analysis, \emph{ Freiburger Intervall - Berichte} 10 (1981), 1 -
63.

\bibitem{Anguelov-Rosinger} R. Anguelov, E.E. Rosinger, Solution of
Nonlinear PDEs by Hausdorff Continuous Functions (to appear).

\bibitem{Baire} R. Baire, Lecons sur les Fonctions Discontinues, Collection
Borel, Paris, 1905.

\bibitem{Bardi} M Bardi, I Capuzzo-Dolcetta, {\it Optimal control and viscosity
solutions of Hamilton-Jacobi-Bellman equations}, Birkh\"{a}user,
Boston, Basel, Berlin, 1997.

\bibitem {Birkhoff}G. Birkhoff, The Role of Order in Computing, in \emph{Reliability in Computing} (ed. R.
Moore) (Academic Press, 1988), 357--378.

\bibitem {Dilworth}R. P. Dilworth, The normal completion of the lattice of
continuous functions, \emph{Trans. Amer. Math. Soc.} {\bf 68}
(1950), 427--438.

\bibitem{Luxemburg} W.A.J. Luxemburg, A.C. Zaanen, Riesz Spaces I, North
Holland, Amsterdam, 1971.

\bibitem {Mack}J. E. Mack and D. G. Johnson, The Dedekind completion of
C(X), \emph{Pacif. J. Math.} {\bf 20} (1967), 231-243.

\bibitem {Markov77}S. Markov, A nonstandard subtraction of intervals, \emph{Serdica} {\bf 3} (1977), 359--370.

\bibitem {Markov2}S. Markov, Calculus for interval functions of a real
variable, \emph{Computing} {\bf 22} (1979), 325--337.

\bibitem {Markov}S. Markov, Extended interval arithmetic involving
infinite intervals, \emph{Mathematica Balkanica} {\bf 6} (1992),
269--304.

\bibitem{Rosinger} M.B. Oberguggenberger, E.E. Rosinger, Solution on
Continuous Nonlinear PDEs through Order Completion, North-Holland,
Amsterdam, London, New York, Tokyo, 1994.

\bibitem {Sendov0}B. Sendov, Approximation of functions by algebraic
polinomials with respect to a metric of Hausdorff type,
\emph{Annals of Sofia University, Mathematics} {\bf 55} (1962),
1-39.

\bibitem {Sendov} B. Sendov, \emph{Hausdorff Approximations}, Kluwer Academic,
Boston, 1990.

\bibitem{Brunn-Minkowski} R. Schneider, \emph{Convex bodies: The
Brunn - Minkowski theory}, Cambridge University Press, 1993.

\end{thebibliography}
\end{document}